\documentclass[reqno, english]{amsart}
\usepackage{etex}
\usepackage{amsmath,amssymb,amsthm,bbm,mathtools,comment}
\usepackage[shortlabels]{enumitem}
\usepackage[pdftex,colorlinks,backref=page,citecolor=blue]{hyperref}
\hypersetup{pdfpagemode=UseNone,pdfstartview={XYZ null null 1.00}}
\usepackage[mathscr]{euscript}
\usepackage[usenames,dvipsnames]{color}
\usepackage{adjustbox,tikz,calc,graphics,babel,standalone}
\usetikzlibrary{shapes.misc,calc,intersections,patterns,decorations.pathreplacing}
\usetikzlibrary{arrows,shapes,positioning}
\usetikzlibrary{decorations.markings}
\usepackage[final]{microtype}
\usepackage[numbers]{natbib}
\usepackage{cmtiup}
\usepackage{amsfonts}
\usepackage{graphicx}
\usepackage{caption}
\usepackage{subcaption}
\usepackage{verbatim}
\usepackage{array}
\usepackage[frame,cmtip,arrow,matrix,line,graph,curve]{xy}
\usepackage{graphpap, color, pstricks}
\usepackage{pifont}
\usepackage[final]{microtype}
\usepackage{cmtiup}
\usepackage{todonotes}
\usetikzlibrary{topaths,calc} 
\tikzstyle{vertex} = [fill,shape=circle,node distance=80pt]
\tikzstyle{edge} = [opacity=0.4,fill opacity=0.0,line cap=round, line join=round, line width=40pt]
\tikzstyle{elabel} =  [fill,shape=circle,node distance=30pt]

\allowdisplaybreaks

\theoremstyle{plain}
\newtheorem{theorem}{Theorem}[section]		
\newtheorem{lemma}[theorem]{Lemma}
\newtheorem{claim}[theorem]{Claim}

\newtheorem{conjecture}[theorem]{Conjecture}
\newtheorem{problem}[theorem]{Problem}

\theoremstyle{remark}

\DeclareMathOperator\ex{ex}

\let\originalleft\left
\let\originalright\right
\renewcommand{\left}{\mathopen{}\mathclose\bgroup\originalleft}
\renewcommand{\right}{\aftergroup\egroup\originalright}

\makeatletter
\def\imod#1{\allowbreak\mkern10mu({\operator@font mod}\,\,#1)}
\makeatother

\title{Two Ramsey problems in blowups of graphs}
\author{Ant\'onio Gir\~ao}
\author{Robert Hancock}
\thanks{
AG: Mathematical Institute, University of Oxford, Oxford OX2 6GG, UK. E-mail: {\tt girao@maths.ox.ac.uk}. {Research supported by Deutsche Forschungsgemeinschaft
(DFG, German Research Foundation) under Germany’s Excellence Strategy EXC-2181/1 - 390900948 (the
14
Heidelberg STRUCTURES Cluster of Excellence) and by EPSRC grant EP/V007327/1}}
\thanks{
RH: Mathematical Institute, University of Oxford, Oxford OX2 6GG, UK. E-mail: {\tt robert.hancock@maths.ox.ac.uk}. Research supported by a Humboldt Research Fellowship at Heidelberg University.}

\begin{document}
\begin{abstract}
 
Given graphs $G$ and $H$, we say $G \stackrel{r}{\to} H$ if every $r$-colouring of the edges of $G$ contains a monochromatic copy of $H$. 
Let $H[t]$ denote the $t$-blowup of $H$. 
The blowup Ramsey number $B(G \stackrel{r}{\to} H;t)$ is the minimum $n$ such that $G[n] \stackrel{r}{\to} H[t]$. 
Fox, Luo and Wigderson refined an upper bound of Souza, showing that, given $G$, $H$ and $r$ such that $G \stackrel{r}{\to} H$, there exist constants $a=a(G,H,r)$ and $b=b(H,r)$ such that for all $t \in \mathbb{N}$, $B(G \stackrel{r}{\to} H;t) \leq ab^t$. 
They conjectured that there exist some graphs $H$ for which the constant $a$ depending on $G$ is necessary. 
We prove this conjecture by showing that the statement is true in the case of $H$ being $3$-chromatically connected, which in particular includes triangles. On the other hand, perhaps surprisingly, we show that for forests $F$, 
there exists an upper bound for $B(G \stackrel{r}{\to} F;t)$ which is independent of $G$.

Second, we show that for any $r,t \in \mathbb{N}$, any sufficiently large $r$-edge coloured complete graph on $n$ vertices with $\Omega(n^{2-1/t})$ edges in each colour contains a member from a certain finite family $\mathcal{F}^r_t$ of $r$-edge coloured complete graphs. 
This answers a conjecture of Bowen, Hansberg, Montejano and M\"uyesser. 
\end{abstract}
\maketitle

\section{Blowup Ramsey numbers}

For graphs $G,H$ and $r \in \mathbb{N}$, we write $G \stackrel{r}{\to} H$ if every $r$-colouring of the edges of $G$ contains a monochromatic copy of $H$.
Given a graph $H$ and $t \in \mathbb{N}$, the \emph{$t$-blowup} of $H$, denoted $H[t]$, is the graph obtained from $H$ by replacing each vertex of $H$ by an independent set of size $t$, and replacing every edge of $H$ by a complete bipartite graph $K_{t,t}$ between the corresponding parts.
We call a copy of $H[t]$ in $G[n]$ \emph{canonical} if it is the $t$-blowup of a copy of $H$ in $G$.
Suppose $G$, $H$ are graphs and $r \in \mathbb{N}$ such that $G \stackrel{r}{\to} H$. 
Then for $t \in \mathbb{N}$, the \emph{blowup Ramsey number} $B(G \stackrel{r}{\to} H;t)$ is the minimum $n$ such that every $r$-colouring of the edges of $G[n]$ contains a monochromatic canonical copy of $H[t]$.

Blowup Ramsey numbers were introduced by Souza in~\cite{Souz}, as a natural generalisation of bipartite Ramsey numbers. 
At the same time, he proved the following upper bound on them.

\begin{theorem}[Souza~\cite{Souz}]
Let $G, H$ be graphs, and let $r \geq 2$ be an integer such that $G \stackrel{r}{\to} H$. There exists a constant $c=c(G,H,r)$ such that for every $t \in \mathbb{N}$, 
\begin{align*}
B(G \stackrel{r}{\to} H;t) \leq c^t.
\end{align*}
\end{theorem}

Souza also proved that the exponential-type bound is necessary, finding a lower bound with an exponential constant which does not depend on $G$. 
This result led him to conjecture that in general the dependence on $G$ is not necessary, that is, for any $H$, the dependence on $G$ in the upper bound can be removed from the function $c$. 

Shortly after, Fox, Luo and Wigderson proved that the exponential constant does not depend on $G$, however, the upper bound still has some dependence on $G$. 
In particular, they proved the following.

\begin{theorem}[Fox--Luo--Wigderson~\cite{FLW}]
Let $G, H$ be graphs, and let $r \geq 2$ be an integer such that $G \stackrel{r}{\to} H$. 
There exist constants $a=a(G,H,r)$ and $b=b(H,r)$ such that for every $t \in \mathbb{N}$, 
\begin{align*}
B(G \stackrel{r}{\to} H;t) \leq ab^t.
\end{align*}
\end{theorem}

However, Fox, Luo and Wigderson conjectured that the dependence on $G$ is necessary for some graphs $H$, contrary to the conjecture of Souza.

\begin{conjecture}[Fox--Luo--Wigderson~\cite{FLW}]
There exists a graph $H$ and integers $r,t \geq 2$ for which the following holds. 
There exist graphs $G_1, G_2, \dots$ such that $G_i \stackrel{r}{\to} H$ for all $i$ and $\sup_i B(G_i \stackrel{r}{\to} H;t) = \infty$. 
\end{conjecture}

In particular Fox, Luo and Wigderson conjectured that this holds with $H$ a triangle and $r=t=2$. 
Our first main result is to confirm this conjecture.

\begin{theorem}~\label{thm:maintriangles}
For every integer $s \geq 2$, there exists a graph $G$ such that $G \stackrel{2}{\to} K_3$ but $G[s] \stackrel{2}{\not\to} K_3[2]$.
\end{theorem}

In fact we can prove the result for a much larger class of graphs. 
A graph is \emph{$3$-chromatically connected} if $|V(G)|\geq 3$ and for all $V \subseteq V(G)$ such that $G[V]$ is bipartite, the graph $G-V$ is connected. 
(In particular, note that a triangle is $3$-chromatically connected.) 
Our theorem extends to all graphs which are $3$-chromatically connected. 
We also note that it makes no difference in the proof to replace $2$ by any $r \geq 2$.
 
\begin{theorem}~\label{thm:main3chrom}
Let $H$ be a $3$-chromatically connected graph and let $r \geq 2$ be an integer. 
For every integer $s \geq 2$, there exists a graph $G$ such that $G \stackrel{r}{\to} H$ but $G[s] \stackrel{r}{\not\to} H[2]$.
\end{theorem}

We prove this result in Section~\ref{subsec:proof1}. 
Our theorem relies on a sparse Ramsey theorem of Ne\v set\v ril and R\"odl~\cite{NR} (see Theorem~\ref{thm:NR}), which in particular, requires the graph in question to be $3$-chromatically connected. 
We are therefore as of yet unable to extend our result to graphs which are not $3$-chromatically connected. 
It would be very interesting to determine precisely which graphs $H$ require some dependence on $G$ in their blowup Ramsey numbers, and which do not.

For an integer $r \geq 2$, we call $G$ \emph{$r$-Ramsey-minimal} for $H$, if $G \stackrel{r}{\to} H$, and for any proper subgraph $G'$ of $G$, we have $G' \stackrel{r}{\not \to} H$. 
A graph $H$ is \emph{$r$-Ramsey-finite} if its class of $r$-Ramsey-minimal graphs is finite. 
Souza~\cite{Souz} proved that $r$-Ramsey-finite graphs do not require the dependence on $G$. 

While any graph containing a cycle is not $r$-Ramsey-finite (as observed by Souza~\cite{Souz}), there are some graphs, for example a star with an odd number of edges, which are $2$-Ramsey-finite~\cite{BEFRS}.
However, Fox, Luo and Wigderson~\cite{FLW} showed that $P_3$, the path with two edges, is not $r$-Ramsey-finite, yet it also does not have any dependence on $G$. 
We are also able to show that the same phenomenon holds for any forest $T$.
\begin{theorem}\label{thm:treesmain}
For any forest $T$ and for any integer $r \geq 2$,
there exists a function $f=f_{T,r}$ such that
for any graph $G$ with $G \stackrel{r}{\to} T$,
we have $B(G \stackrel{r}{\to} T;t) \leq f_{T,r}(t)$ for all $t \in \mathbb{N}$.
\end{theorem}
Note that the blowup Ramsey number of a non-connected graph is clearly at most the sum of the blowup Ramsey numbers for each of its components, and therefore it suffices to prove the above result for trees.
We prove this result in Section~\ref{subsec:prooftrees}.

\subsection{Proof of Theorems~\ref{thm:maintriangles} and~\ref{thm:main3chrom}}\label{subsec:proof1}  \hfill\\
Given two graphs $H$ and $G$, we define $H(G)$, the \emph{hypergraph of copies of $H$ in $G$}, to be the $|E(H)|$-uniform hypergraph with vertex set given by the edges $E(G)$ of $G$, and edge set consisting of each set of edges of $G$ which correspond to a copy of $H$ in $G$.  

We define a \emph{cycle} in a hypergraph $H$ to be a set of edges $e_1,\dots,e_k$ such that there are distinct vertices $v_1,\dots,v_k$ for which $v_i \in e_i \cap e_{i+1}$ for each $i \in [k-1]$ and $v_k \in e_1 \cap e_k$.
In particular, note that if two edges share at least two vertices, then these two edges form a cycle of length $2$. 
We define the \emph{girth} of a hypergraph $H$ to be the number of edges in a shortest cycle of $H$ and denote it by $g(H)$.
Note therefore that a hypergraph with girth at least three is what is often referred to in literature as a \emph{linear} hypergraph, and all its cycles are \emph{loose} cycles, i.e. each pair of consecutive edges along the cycle intersect in precisely one vertex.

The graph $G$ which we use in the proof will be obtained using the following sparse Ramsey theorem of Ne\v set\v ril and R\"odl.

\begin{theorem}[Ne\v set\v ril--R\"odl, Theorem 4.2~\cite{NR}] \label{thm:NR}
Let $r,s \in \mathbb{N}$ and let $H$ be a $3$-chromatically connected graph. 
Then there exists a graph $G$ with the following properties:
\begin{enumerate}
\item $G \stackrel{r}{\to} H$;
\item $g(H(G)) > s$.
\end{enumerate}
\end{theorem}

The only property of a $3$-chromatically connected graph we require in our proof is the existence of the above $G$, together with having minimum degree at least $2$ (which is required only at stage 1, step (i)(1) of the proof).

So for the rest of the proof, to make it clearer and easier to read, we will take $r=2$ and $H=K_3$, noting that it easily extends to general $r \geq 2$ and $H$ which are $3$-chromatically connected.

Fix $s \geq 2$ and choose $G$ according to Theorem~\ref{thm:NR} so that $g(H(G))>2s+2$. 
Note that we can choose $G$ so that it is ($2$-)Ramsey-minimal, which we recall means that for any edge $e$, we have $G-e \stackrel{2}{\not \to} K_3$.  
In order to complete the proof, it suffices to find a red/blue colouring of $G[s]$ which does not contain any monochromatic canonical copies of $K_3[2]$.

We first give an outline of the proof.
We start with a colouring of $E(G)$ which only contains one monochromatic copy of $K_3$, and transfer this colouring to $E(G[s])$ in the obvious way, colouring every edge in the $K_{s,s}$ blowup of $uw$ by the colour of $uw$.
The idea is to step by step recolour edges of $G[s]$ so that in each step all monochromatic copies of $K_3[2]$ created by the previous step are destroyed, while in each edge blowup, less edges are recoloured compared to the previous step. 
In stage 1, we will look at colourings of $G$.
Since $G$ is Ramsey-minimal, we can select an edge $e \in E(G)$, for which there exists a red/blue-colouring of $E(G-e)$ such that there are no monochromatic triangles. 
Now, we extend the red/blue-colouring of $E(G-e)$ to a red/blue-colouring of $E(G)$ by colouring $e$ red. 
Since $G \stackrel{2}{\to} K_3$, there must exist some red triangles; note that all of these red triangles contain $e$, and further there are no blue triangles. 
Further note that these red triangles intersect at $e$ and are otherwise edge-disjoint from each other, since $H(G)$ contains no cycle of length $2$. 
(Note that actually two triangles can only intersect in one edge in the first place, but this argument ensures that when we replace triangles by general $H$, it is still the case that these copies of $H$ are edge-disjoint from each other except for at $e$.)
Therefore we can select precisely one edge per red triangle to recolour blue. 
We obtain a colouring which no longer contains red triangles, although it may contain new blue triangles, each of which must contain one of the recoloured edges. 
Note that again since $H(G)$ does not contain short cycles, the set of blue triangles each share at most one edge with any other triangle (precisely the edge which was just recoloured from red to blue) and are otherwise entirely edge-disjoint from all other monochromatic triangles created previously. 
We do this recolouring step $s$ times. 
In stage 2, we will transfer these edge colourings of $G$ to $G[s]$ so that in each step, whenever we recolour edges of $G$, we recolour proportionally less of the edges in the blowup of these edges of $G$ compared to the previous step, while still destroying all monochromatic copies of $K_3[2]$ created during the previous step.
Eventually we get to a step where we do not have to recolour anything at all in order to destroy all monochromatic copies of $K_3[2]$.  

We now make this recolouring procedure precise. 

{\bf Stage 1: the colourings in $G$.}

{\bf Step 0:}
Fix an edge $e \in E(G)$ and let $c_0$ be a red/blue colouring of $E(G)$ such that there are no monochromatic triangles in $E(G-e)$ and $e$ is red. 
Further, fix a vertex $v$ incident to $e$, and set $E_0=\{e\}$.
Let $T_0$ be the set of triangles in $G$ which are red in $c_0$. 

Note that all triangles in $T_0$ contain $e$, and further except for $e$, these triangles are edge-disjoint. This follows from $g(H(G))>2s+2$. 

We now repeat the following step $s$ times, for each $i=1,\dots,s$.

{\bf Step $i$:}
\begin{enumerate}
\item For each $T \in T_{i-1}$, change the colour of the edge of $T$ which is not in $E_{i-1}$ and is incident to $v$. (Note that for general $H$ and $r$, there are many such choices of edge or colour, but any choice works.)
\item Call the new colouring $c_i$.
\item Call the set of edges which were just recoloured $E_i$.
\item Let $T_{i}$ be the set of monochromatic triangles in $c_i$.
\end{enumerate}

Note that after each step, the following holds.
In $c_i$, there are only monochromatic triangles of one colour; in particular, a colour which did not occur in the previous colouring $c_{i-1}$.
Every edge in $E_i$ contains the fixed vertex $v$.
All triangles in $T_i$ contain exactly one edge from $E_{i}$ (see the edge-disjointedness argument below) and all of these triangles contain the vertex $v$.


Crucially, all triangles in $\cup_{j=0}^{i} T_j$ are edge-disjoint from each other, except for the following cases:
\begin{itemize} 
\item A triangle $A_i \in T_i$ shares precisely one edge $f$ with one triangle $A_{i-1} \in T_{i-1}$, where $f \in E_i$;
\item A triangle $A_i \in T_i$ may share precisely one edge $f$ with other triangles in $T_i$; if so, then $f \in E_i$.
\end{itemize}

This edge-disjointedness follows from $g(H(G))>2s+2$; 
suppose for a contradiction that there exists triangles $A_i \in T_i$ and $B_j \in T_j$, $j \leq i$, which violate the above. 
If they share two edges then this immediately contradicts $g(H(G))>2s+2$ since they would form a cycle of length $2$ in $H(G)$. 
Otherwise, by design, for each $0 \leq k <i$ we can find triangles $A_k \in T_k$ such that $A_k \cap A_{k+1} = e_{k+1} \in E_{k+1}$, and for each $0 \leq k < j$ we can find triangles $B_k \in T_k$ such that $B_k \cap B_{k+1} = f_{k+1} \in E_{k+1}$. 
Let $\ell$ be the maximum index for which $A_{\ell}=B_{\ell}$. 
If $\ell$ exists then the triangles $A_{\ell}, A_{\ell+1}, \dots, A_{i}, B_{j}, B_{j-1}, \dots, B_{\ell+1}, A_\ell$ form a cycle in $H(G)$. 
If not, then the triangles $A_0, \dots, A_{i}, B_j, B_{j-1}, \dots, B_0, A_0$ form a cycle in $H(G)$. 
In both cases we have a cycle of length at most $2s+2$, a contradiction to $g(H(G))>2s+2$.
We get a similar contradiction if we assumed that a triangle in $T_i$ contains more than one edge from $E_i$.

{\bf Stage 2: transferring the colourings to $G[s]$.}

{\bf Step 0:}
First we transfer the colouring $c_0$ of $E(G)$ to a colouring $c_0'$ of $E(G[s])$ in the obvious way; 
that is, colour every edge in the $K_{s,s}$ blowup of $uw$ by $c_0(uw)$ for all edges $uw \in E(G)$.

For the fixed vertex $v$, label the vertices in its blowup by $v_1,\dots,v_s$.
We now repeat the following step $s$ times, for each $i=1,\dots,s$.

{\bf Step i:}
\begin{enumerate}
\item For each $uv \in E_i$ (recall $v$ is fixed and $c_{i-1}(uv) \not= c_i(uv)$), recolour edges from $\{v_{i+1},\dots,v_{s}\}$ to all vertices in the blowup of $u$. (Do not change the colour of edges involving $\{v_1,\dots,v_{i}\}$.)
\item Call the new colouring $c'_i$.
\end{enumerate}

\begin{figure}[h]
\begin{center}
\begin{tikzpicture}[scale=0.9]
\filldraw[color=black!100, fill=red!0, very thick](3,4.2) ellipse (2.5 and 1.5);
\node[vertex] at (1,4) (v1) {};
\node[vertex] at (1.75,4) (v2) {};
\node[vertex] at (2.5,4) (v3) {};
\node[vertex, draw, fill=white, inner sep=10pt] at (4.2,4) (v4) {};
\node at (1,4.4) {$v_1$};
\node at (1.75,4.4) {$v_2$};
\node at (2.5,4.4) {$v_3$};
\node at (4.2,4.7) {$\{v_4,\dots,v_s\}$};
\node at (5.7,4.4) {$v$};

\node[vertex, draw, fill=white, inner sep=10pt] at (0,1) (x) {};
\node[vertex, draw, fill=white, inner sep=10pt] at (2,0) (y) {};
\node[vertex, draw, fill=white, inner sep=10pt] at (4,0) (u) {};
\node[vertex, draw, fill=white, inner sep=10pt] at (6,1) (w) {};

\node at (0,1) {$x$};
\node at (2,0) {$y$};
\node at (4,0) {$u$};
\node at (6,1) {$w$};

\draw[very thick,red] (v1) -- (x);
\draw[very thick,red] (v1) -- (y);
\draw[very thick,blue] (v1) -- (u);
\draw[very thick,red] (v1) -- (w);

\draw[very thick,red] (v2) -- (x);
\draw[very thick,blue] (v2) -- (y);
\draw[very thick,blue] (v2) -- (u);
\draw[very thick,red] (v2) -- (w);

\draw[very thick,red] (v3) -- (x);
\draw[very thick,blue] (v3) -- (y);
\draw[very thick,red] (v3) -- (u);
\draw[very thick,red] (v3) -- (w);

\draw[very thick,red] (v4) -- (x);
\draw[very thick,blue] (v4) -- (y);
\draw[very thick,red] (v4) -- (u);
\draw[very thick,red] (v4) -- (w);

\draw[very thick,red] (x) -- (y);
\draw[very thick,blue] (y) -- (u);
\draw[very thick,red] (u) -- (w);

\draw[->,very thick] (6,3) -- (8,3);

\filldraw[color=black!100, fill=red!0, very thick](11,4.2) ellipse (2.5 and 1.5);
\node[vertex] at (9,4) (v1b) {};
\node[vertex] at (9.75,4) (v2b) {};
\node[vertex] at (10.5,4) (v3b) {};
\node[vertex, draw, fill=white, inner sep=10pt] at (12.2,4) (v4b) {};
\node at (9,4.4) {$v_1$};
\node at (9.75,4.4) {$v_2$};
\node at (10.5,4.4) {$v_3$};
\node at (12.2,4.7) {$\{v_4,\dots,v_s\}$};
\node at (13.7,4.4) {$v$};

\node[vertex, draw, fill=white, inner sep=10pt] at (8,1) (xb) {};
\node[vertex, draw, fill=white, inner sep=10pt] at (10,0) (yb) {};
\node[vertex, draw, fill=white, inner sep=10pt] at (12,0) (ub) {};
\node[vertex, draw, fill=white, inner sep=10pt] at (14,1) (wb) {};

\node at (8,1) {$x$};
\node at (10,0) {$y$};
\node at (12,0) {$u$};
\node at (14,1) {$w$};

\draw[very thick,red] (v1b) -- (xb);
\draw[very thick,red] (v1b) -- (yb);
\draw[very thick,blue] (v1b) -- (ub);
\draw[very thick,red] (v1b) -- (wb);

\draw[very thick,red] (v2b) -- (xb);
\draw[very thick,blue] (v2b) -- (yb);
\draw[very thick,blue] (v2b) -- (ub);
\draw[very thick,red] (v2b) -- (wb);

\draw[very thick,red] (v3b) -- (xb);
\draw[very thick,blue] (v3b) -- (yb);
\draw[very thick,red] (v3b) -- (ub);
\draw[very thick,red] (v3b) -- (wb);

\draw[very thick,red] (v4b) -- (xb);
\draw[very thick,blue] (v4b) -- (yb);
\draw[very thick,red] (v4b) -- (ub);
\draw[very thick,blue] (v4b) -- (wb);

\draw[very thick,red] (xb) -- (yb);
\draw[very thick,blue] (yb) -- (ub);
\draw[very thick,red] (ub) -- (wb);

\end{tikzpicture}
\end{center}
\caption{An example of step 3 of stage 2. 
The triangle $uvw$ is the monochromatic $T \in T_{2}$, with $uv \in E_2$ and $vw \in E_3$. 
Originally (on the left) all edges between vertices in the blowups of $v$ and $w$ are red. 
In step 3, we recolour all edges between $\{v_4,\dots,v_s\}$ and vertices in the blowup of $w$ from red to blue (on the right). It is not too hard to see that there are no monochromatic canonical copies of $K_3[2]$ with parts in $u,v,w$.  
Note we can see evidence of the previous recolourings too;
in step 2, all edges between $\{v_3,\dots,v_s\}$ and vertices in the blowup of $u$ were coloured red;
in step 1, all edges between $\{v_2,\dots,v_s\}$ and vertices in the blowup of $y$ were coloured blue.}
\label{fig:recol}
\end{figure}
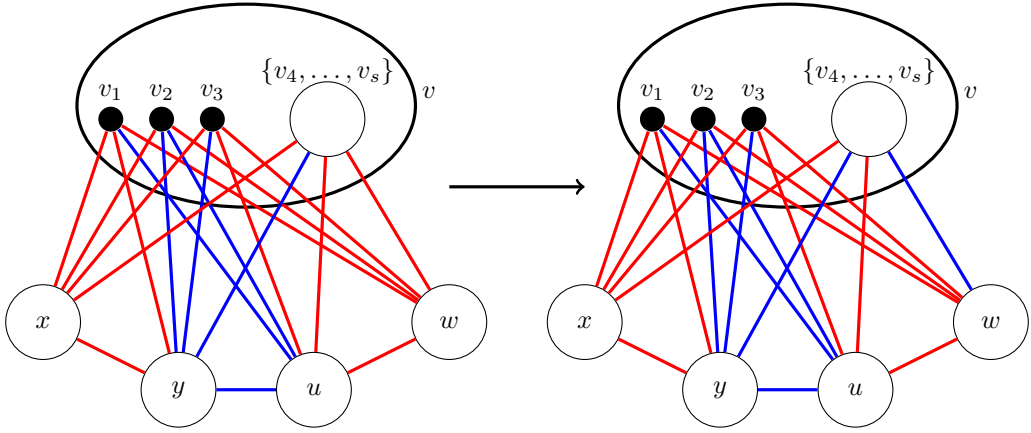

See Figure~\ref{fig:recol} for an example of step 3.

We now aim to prove that $c'_{s}$ is a colouring of $G[s]$ which contains no monochromatic canonical copies of $K_3[2]$. This will follow quickly from the following claim.

\begin{claim}\label{claim:recol}
In step $i$, for $i \in [s]$, all monochromatic canonical copies of $K_3[2]$ contained in $\{T[s] : T \in T_{i-1}\}$ are destroyed.
\end{claim}
 
Indeed, recall that in step $i$ of stage 1, all monochromatic copies of $K_3$ are collected in the set $T_i$. 
Therefore, by Claim~\ref{claim:recol},  
the only possible monochromatic copies of $K_3[2]$ after step $i$ of stage 2 must contain edges which were just recoloured within this step. 
Since in step $s$ of stage 2 no edges are recoloured, the colouring $c'_{s}$ contains no monochromatic copies of $K_3[2]$ as desired.
It remains to prove the claim.

\begin{proof}[Proof of Claim~\ref{claim:recol}]
Since all $T \in T_{i-1}$ contain $v$, we must use $2$ vertices from $\{v_1,\dots,v_s\}$ in any monochromatic $K_3[2]$. 
Since all edges involving $\{v_{i+1},\dots,v_{s}\}$ are recoloured, all monochromatic copies of $K_3[2]$ using these vertices are destroyed, and so we must use $2$ vertices from $\{v_1,\dots,v_{i}\}$. 
In step 1, we are done, since there are not enough vertices.

Now consider any step $i$ with $i \geq 2$.
By construction, all triangles $T_{i-1}$ are of the form $uvw$, with $vw \in E_i$ and $uv \in E_{i-1}$. 
Since $uv \in E_{i-1}$, in step $i-1$, the colour of the edges $\{v_1,\dots,v_{i-1}\}$ to vertices in the blowup of $u$ were not changed, and hence are a different colour to the edges from $\{v_1,\dots,v_i\}$ to vertices in the blowup of $w$.
Therefore $v_i$ is the only vertex with the same colour of edges to vertices in the blowups of $u$ and $w$, and hence we cannot find a monochromatic $K_3[2]$ using the vertices $u$, $v$ and $w$. 
\end{proof}

\subsection{Proof of Theorem~\ref{thm:treesmain}}\label{subsec:prooftrees} \hfill\\
We will prove a stronger result, Theorem~\ref{thm:treesstronger}, from which our main result Theorem~\ref{thm:treesmain} will immediately follow. Note that as with the proof of Theorem~\ref{thm:main3chrom}, we will assume $r=2$ and note that the proof easily extends to $r \geq 3$.

First we need to introduce some notation.
Let $G$ be a graph. 
We say $c$ is a \textit{partial $2$-colouring} of $G$ if there exists a subgraph 
$G'\subset G$ such that $c: E(G')\rightarrow \{\text{red},\text{blue}\}$, 
i.e. not necessarily all edges receive colours, but those which do are one of two colours.

Given a tree $T$, a graph $G$ and a partial $2$-colouring $c$ of $G$, 
we say that $H\subset G$ is a \textit{possible red (blue) copy} of $T$ if 
$H$ is isomorphic to $T$ and 
all edges of $H$ coloured by $c$ are red (blue).
(Note that not all edges of $H$ need to be coloured.) 
We say $H$ is a \textit{possible monochromatic copy} of $T$ if $H$ is either a possible red or possible blue copy of $T$. 

Let $G$ be a graph, as usual we denote by $G[m]$ the $m$-blowup of $G$ and 
more generally,
given a function $m: V(G) \to \mathbb{N}$,
we denote by $G[\{m(x)\}_{x\in V(G)}]$ the blowup of $G$ where 
each vertex $x$ is blown up by $m(x)$ vertices.
For every $x\in V(G)$, we denote by $\overline{x}$ the blowup of $x$ in $G[\{m(x)\}_{x\in V(G)}]$.
We write $(\overline{x},\overline{y})$ to denote the complete bipartite graph between $\overline{x}$ and $\overline{y}$ in $G[\{m(x)\}_{x\in V(G)}]$.
  
Let $G$ be a graph, $c$ a partial $2$-colouring of $G$ and 
let $f:\mathbb{N}\rightarrow \mathbb{N}$ be a non-decreasing function such that $\lim_{n\rightarrow \infty} f(n)=\infty$.
(We introduce this function so that we get a cleaner argument which avoids tracing all of the constants.)
We say a colouring $c': E(G[\{m(x)\}_{x\in V(G)}])\rightarrow \{\text{red}, \text{blue}\}$ is 
$f$-\textit{coherent with} $c$ if for every edge $e=xy \in E(G)$ if $c(e)$ is red (blue) 
then the complete bipartite graph $(\overline{x},\overline{y})$ does not contain a blue (red) $K_{f(m(x)),f(m(y))}$ in $c'$. 
Note that we will always choose $f$ so that $K_{a,b} \to K_{f(a),f(b)}$ for all $a,b \in \mathbb{N}$, so in particular, if the colour of an edge $xy$ is red (blue), then, since $(\overline{x},\overline{y})$ does not contain a blue (red) $K_{f(m(x)),f(m(y))}$, it must contain a red (blue) $K_{f(m(x)),f(m(y))}$.

We now prove a technical lemma. 
%

\begin{lemma}\label{lem:tree2}
Let $T$ be a tree, let $T'\subset T$ be a subtree, and let $xy\in E(T)\setminus E(T')$ where $x\in V(T')$. 
Let $G$ be a graph, and let $c$ be a partial $2$-colouring of $G$.
Let $z\in V(G)$, and let $A=\{T'_1,\ldots T'_k\}$ be a set of distinct copies of $T'$ in $G$ such that $z$ plays the role of $x$ in all of these copies of $T'$,
and for every $i\in [k]$ there is a possible monochromatic copy of $T$ in $G$, 
say $T_i$, such that $T'_i\subset T_i$. 
If $k$ is sufficiently large, then there exist $i,j$ with $i<j$ such that $T_i$ and $T_j$ are possible monochromatic copies of $T$, and $T_j \setminus T'_j$ is vertex disjoint from $T'_i$.  
\end{lemma}

\begin{proof}
First we may pass to a subset $A'=\{T'_1,\ldots, T'_{k_0}\}$ of $A$ of at least half the size where 
without loss of generality we may assume all copies of $T'_i$ lie in a possible red copy of $T$. 
If $T_j\setminus T'_j$ is vertex disjoint from $T'_1$ 
for some $1 < j \leq k_0$, then we are done.
Therefore, we may assume $T_j\setminus T'_j$ has non-empty intersection with $T'_1$ for all $1 < j \leq k_0$. 
Let $x_1$ be a vertex of $T'_1$ which appears most often (at least $\frac{k_0}{v(T')}$ times) in the intersection of $T'_1$ with $T_j\setminus T'_j$, $1 < j \leq k_0$.
We may pass to a subset of $A'$, say without loss of generality $A_1=\{T'_2,\ldots , T'_{k_1}\}$, such that 
for every $2 \leq \ell \leq k_1$, $T_{\ell}\setminus T'_{\ell}$ contains the vertex $x_1$. 
Now if $T_j\setminus T'_j$ is vertex disjoint from $T'_2$ for some $2<j\leq k_1$, then we are done.  
Hence, again, we may pass to a subset of $A_1$, say without loss of generality $A_2=\{T'_3,\ldots , T'_{k_2}\}$, such that for every $3 \leq \ell \leq k_2$, $T_{\ell}\setminus T'_{\ell}$ contains the vertex $x_2$, 
which was a vertex of $T'_2$ appearing most often in the intersection of $T'_2$ with $T_j\setminus T'_j$, $2 < j \leq k_1$.
(Note that $x_2 \not= x_1$ since $x_1 \in T_2 \setminus T'_2$.)
We keep doing this process, forming further subsets $A_3, A_4, \dots$, and finding further vertices $x_3, x_4, \dots$ which lie in each $T_{\ell}\setminus T'_{\ell}$.
If $k$ is sufficiently large (in particular $k \geq 2v(T')^{v(T)-v(T')}$ suffices)
then this process must stop, and we have the conclusion.
\end{proof}

We are now ready to prove the following result which easily implies Theorem~\ref{thm:treesmain}. It turns out to be more convenient to allow each vertex to be blown up a different number of times. 
\begin{theorem}\label{thm:treesstronger}
For every tree $T$, and 
for any subtree $T' \subset T$, 
there exists a non-decreasing function $f$ with $\lim_{n\rightarrow \infty} f(n)=\infty$ 
such that the following holds. 
Let $m: V(G) \to \mathbb{N}$,
let $G \stackrel{2}{\to} T$ 
and let $c$ be a partial $2$-colouring of $G$. 
Then, for any $2$-colouring $c'$ of $E(G[\{m(x)\}_{x\in V(G)}])$ which is $f$-coherent with $c$,
there exists $H,H'\subset G$ 
such that $H'\subset H$, 
$H'$ is a copy of $T'$, 
and $H$ is a possible red (or blue) copy of $T$ in $G$ using $c$,
and there exists a canonical copy of $H'[\{f(m(x))\}_{x\in V(H')}]$ monochromatic in red (or blue) in $G[\{m(x)\}_{x\in V(G)}]$ using $c'$.
\end{theorem}

Note that Theorem~\ref{thm:treesmain} for $r=2$ and trees easily follows by applying this theorem with $T'=T$, letting $c$ be empty and by choosing $m$ in hindsight so that for all $x \in V(G)$, $m(x)=m$ and $f(m) \geq t$. 

Before we proceed with the proof, we give a brief high-level description. 
We will proceed by induction on $|E(T')|$. 
The idea is that we take the canonical monochromatic blowup of $T''$, where $T''=T' \setminus y$ for some leaf $y$ of $T'$, which exists by induction, and try to extend it to a canonical monochromatic blowup of $T'$. 
We look at all possible ways of extending $T''$ to some $T'$ within some copy of $T$ in $G$, noting that if we fail for a particular pair $T'$ and $T''$, it must be because the large monochromatic bipartite graph in the blowup of the edge $e$ which we tried to extend along (which must exist by Ramsey's theorem) is monochromatic in the wrong colour. 
We then extend the partial $2$-colouring of $G$ by marking $e$ by this wrong colour and we shrink the size of the blowup of the vertex we were trying to extend from. 
Then we start the inductive step again. 
Shrinking the size of the blowup ensures that we do no select the same $T''$ to try to extend from as previously.
Each time we fail, we colour the edges for which we failed to extend along and shrink the size of the blowups and then try again.
Lemma~\ref{lem:tree2} is used to guarantee that this process must stop, in particular for each vertex the size of its blowup can only be shrunk a finite number of times. Therefore we must eventually succeed and obtain a monochromatic bipartite graph of the correct colour. 
Note that Lemma~\ref{lem:tree2} only holds for trees.

\begin{proof}
First fix the larger tree $T$. We prove by induction on $|E(T')|$. 
For the base case $T'$ is a single edge. 
Let $f$ be a non-decreasing function $f$ such that $\lim_{n\rightarrow \infty} f(n)=\infty$ and which grows sufficiently slowly so that
$K_{a,b} \stackrel{2}{\to} K_{f(a),f(b)}$ for all $a,b \in \mathbb{N}$. 
Let $c,c'$ be as in the statement. 
Since $G \stackrel{2}{\to} T$ there exists a possible monochromatic copy of $T$ in $G$.
If it contains a red (blue) edge $xy$, then since $c'$ is $f$-coherent with $c$, 
then by Ramsey's Theorem  (and by the definition of $f$), $(\overline{x},\overline{y})$ contains a canonical red (blue) $T'[\{ f(m(x)),f(m(y)) \}]$ as required.
If not, then $T$ does not have any colours and is a possible red copy and a possible blue copy.
Then applying Ramsey's Theorem to any edge $xy$ of $T$, one can find a red or blue canonical $T'[\{ f(m(x)),f(m(y)) \}]$ (it does not matter which colour), as required.

Now suppose that the result holds for $T''$, where $T''=T' \setminus y$ and $y$ is a leaf of $T'$ and we wish to prove the result for $T'$.
Suppose $y$ is adjacent to $x$ in $T'$. 
We need to show that the function $f$ exists. 
By induction with $T,T''$ there is a function $g$ satisfying the properties in the statement. 
Let $f=g^{k}$ i.e. the composition of $g$ $k$ times, where $k$ is the number outputted by Lemma~\ref{lem:tree2} using $T$ and $T''$. 

By induction we may find a canonical $H''[\{g(m(v))\}_{v\in V(H'')}]$ monochromatic in $c'$  
where $H''\subset H$ and $H$ is a possible red (blue) copy of $T$ and $H''$ is a copy of $T''$. 
(Note that $c'$ is $g$-coherent with $c$ because $f \leq g$, which follows from the fact that $g$ is a non-decreasing function with $g(a) \leq a$ for all $a \in \mathbb{N}$, and $f$ is the composition of $g$ $k$ times. We have $g(a) \leq a$ since this is true for the function from the base case.) 
Let $\{H_1,\ldots, H_t\}$ be all possible red (blue) copies of $T$, where for each $i \in [t]$, $H_i$ contains $H'_i$, a copy of $T'$, where $T''=T' \setminus y$ and $y$ is a leaf of $T'$, and each $H'_i$ contains $H''$.
(In particular, we allow any way of embedding $T'$ into $T$, that is, we do not need each $H'_i$ to play the same role as a subgraph of $T$, we only require that $H''$ can be extended to a copy of $T'$ within each $H_i$.)
For each $i \in [t]$, let $y_i\in V(H_i)$ be the vertex which corresponds to the vertex $y$. 
Let $H''(x)$ be the subset of $\overline{x}$ of size $g(m(x))$ which lies in the canonical monochromatic $H''[\{g(m(v))\}_{v\in V(H'')}]$ in $c'$.  
If we can find a red (blue) copy of $K_{g(g(m(x))),g(m(y_i))}$ between $H''(x)$ and a subset of size $g(m(y_i))$ in $\overline{y_i}$ (for some $i$) we are done;
this red (blue) complete bipartite subgraph of $(H''(x),\overline{y_i})$ together with the canonical monochromatic $H''[\{g(m(v))\}_{v\in V(H'')}]$ forms a canonical monochromatic $H'[f(m(v))_{v\in V(H')}]$ 
where $H'' \subset H' \subset H_i=H$, $H$ is a possible red (blue) copy of $T$ and $H'$ is a copy of $T'$, as required.
 
So suppose this is not the case. In particular this means that no edge $xy_i$ was coloured red (blue) in $c$; 
as otherwise, since $c'$ is $g$-coherent with $c$, by applying Ramsey's Theorem to $(H''(x),\overline{y_i})$ we would have been able to find a red (blue) $K_{g(g(m(x))),g(m(y_i))}$ in $(H''(x),\overline{y_i})$. These edges also cannot be blue (red) in $c$ since they are supposed to form part of a red (blue) possible copy of $T$. 
Therefore these edges must be uncoloured.


Hence, we may shrink $\overline{x}$ such that only the vertices in $H''(x)$ remain and we denote the new blowup of $G$ by $G^*_1$. 
We also update $c$ by colouring all edges $xy_i$ blue (red) and denote it by $c_1$. 
Clearly, the colouring $c'$ of $G^*_1$ is still $g$-coherent with $c_1$, by construction. 
Also note that now $H''$ cannot be a subgraph of a possible monochromatic copy of $H$. 
Further, without replacing $\overline{x}$ by $H''(x)$, it may not be true that $c'$ is $g$-coherent with $c_1$,
since there may be a large red (blue) bipartite graph between $\overline{x}$ and $\overline{y_i}$ for some $i$. 
We only know that no such bipartite graph exists between $H''(x)$ and $\overline{y_i}$. 
Finally, we define a new function $m_1:V(G) \to \mathbb{N}$ by defining $m_1(x)=g(m(x))$ and $m_1(y)=m(y)$ for all $y \in V(G) \setminus \{x\}$.

Now, we keep applying induction, using the same $T''$, each time shrinking the blowup of one vertex of $G$, and adding new colours to $c$ and updating the function $m$. 
(For simplicity we keep the same names $G$, $c$ and $m$.)
Note that after shrinking a blowup of a vertex $\overline{x}$ for the $i$-th time, it goes from size $g^{i-1}(m(x))$ to $g^{i}(m(x))$.

The essential fact is that each blowup of a vertex will not be shrunk more than $k$ times because of Lemma~\ref{lem:tree2}.
To see this, observe the following: suppose that $\{H''_1, H''_2, \dots,H''_{k}\}$ are each canonical copies of $T''$ found by induction each using the same vertex $x$. 
Let $H_i$ be the possible monochromatic copy of $T$ which $H''_i$ lies in at the point when $H''_i$ was first found by induction for each $i \in [k]$.
Note that after $c$ is updated, $H_i$ is no longer a monochromatic copy.
However, crucially, in the initial colouring $c$, $H_i$ is a possible monochromatic copy which $H''_i$ lies in, 
since by removing colours of edges, we do not destroy any possible monochromatic copies of $T$.
For each $i \in [k]$, let $y_i \notin V(H''_i)$ be a neighbour of $x$ in $H_i$ playing the same role in $T$.
Let $T_{y_i}$ be the subtree of $H_i$ rooted at $y_i$. 
Let $J_i$ be the tree obtained from $H_i$ by removing $T_{y_i}$ and the edge $xy_i$.
Now apply Lemma~\ref{lem:tree2} with $z=x$, $T'_i=J_i$, $T_i=H_i$ and the original colouring $c$.
Since $k$ is large enough, the lemma outputs $i<j$ such that $H_i$ and $H_j$ are both possible red (blue) copies of $T$ and $J_i$ is vertex disjoint from $H_j \setminus J_j$.
Now observe that $H_j \setminus J_j=T_{y_j}$, and therefore the tree $J_i \cup T_{y_j} \cup \{xy_j\}$ is isomorphic to $T$,
and further, it is a red (blue) possible copy of $T$ which contains $H''_i$.
As the process did not stop after finding $H''_i$, we must therefore have updated $c$ by colouring $xy_j$ blue (red).
Now when we apply induction and find $H''_j$, since $H_j$ contains $xy_j$, $H_j$ must be a possible blue (red) copy of $T$, a contradiction, since we know it should be a possible red (blue) copy.

Since each blowup of a vertex is shrunk at most $k$ times and each application of induction results in the shrinking of one blowup of a vertex, this process must stop and therefore we must have found the desired monochromatic canonical blowup of $T'$. 
\end{proof}

\section{Unavoidable colour patterns}
Ramsey's Theorem states that for any $r,t \in \mathbb{N}$ and $n$ sufficiently large, 
any $r$-colouring of the edges of $K_n$ contains a monochromatic copy of a clique on $t$ vertices. 
In other words, Ramsey's Theorem states that a monochromatic colouring is an unavoidable colour pattern of a clique. 
It is therefore natural to ask whether there are other colour patterns of a clique which are unavoidable. 
Of course, the answer is negative since the colouring of the large clique might itself be monochromatic. 
One may then restrict the attention to a proper subset of the collection of all $r$-edge colourings of a clique and ask whether within that class there are other colour patterns which are unavoidable. 

First consider the case of $r=2$. 
We say that a $2$-edge coloured $K_{2t}$ is \emph{$(2,t)$-unavoidable} if one of the colours induces either a clique on $t$ vertices or two vertex disjoint cliques on $t$ vertices.
Bollob\'as asked whether for every $\varepsilon>0$ and $t\in \mathbb{N}$, there is $n(\varepsilon,t)$ such that for every $n\geq n(\varepsilon,t)$, 
every $2$-coloured clique on $n$ vertices, where each colour appears in at least $\varepsilon \binom{n}{2}$ edges 
must necessarily contain a $(2,t)$-unavoidable coloured $K_{2t}$. 
This was confirmed by Cutler and Mont\'agh~\cite{CM}. 
Note that the four colour patterns in the family of $(2,t)$-unavoidable graphs are the only such family, since any other colour pattern can be avoided by making the colouring of $K_n$ of this form, and as remarked by Cutler and Mont\'agh, it is minimal, since if one of the four colour patterns in the family is left out, then there is a colouring of $K_n$ which avoids the other three.

Later, Fox and Sudakov~\cite{FS} proved $n(\varepsilon,t)=\varepsilon^{-O(t)}$, which is tight up to the implied constant. 
Instead of imposing both colours to appear in a positive proportion of the edges, we could ask the following.
\begin{problem}
Let $t\in \mathbb{N}$, and let $n_0(t)$ be a sufficiently large integer. Given  $n\geq n_0$, how many edges $m(t,n)$ of each colour must exist in a $2$-edge colouring of $K_n$ to guarantee the existence of a $(2,t)$-unavoidable pattern?
\end{problem}
Note that the Cutler and Mont\'agh result implies that $m(t,n)=o(n^2)$, for fixed $t$ as $n\rightarrow \infty$. 
Quite recently, Caro, Hansberg and Montejano~\cite{CHM1} proved that for every $t \in \mathbb{N}$, there exists a $\delta=\delta(t)$ such that $m(t,n)=O(n^{2-\delta})$, 
following which the first author and Narayanan~\cite{GN} showed that $\delta=1/t$, by proving that for every $t\in \mathbb{N}$, there exists a constant $C(t)$ such that $m(t,n)\leq C(t)n^{2-1/t}$. 
It is not hard to see that this bound is tight up to the constant $C(t)$,
conditional on a famous conjecture of K\H{o}v\'ari-S\'os-Tur\'an~\cite{KST},
which asserts that $\ex(n,K_{a,b})=\Omega(n^{2-1/a})$ for all $b\geq a \geq 2$.

Very recently, Caro, Hansberg and Montejano have obtained various generalisations, as well as cases of unconditional sharpness, of this result, see~\cite{CHM2} for the precise results.
Here, we instead consider the generalisation of the above problem to $r$ colours.
The concept of a $(2,t)$-unavoidable colour pattern was first generalised to $r$ colours by Bowen, Lamaison and M\"uyesser in~\cite{BLM}. 
Let $G$ be a vertex and edge coloured complete graph with $r$ colours. 
Moreover, suppose that all $r$ colours are present and $G$ is minimal with this property, 
i.e. no proper induced subgraph of $G$ spans all $r$ colours. 
Then, we say $G$ is an \emph{$r$-minimal graph}.
The \emph{$t$-blowup} of a vertex and edge coloured graph $G$ on $k$ vertices is obtained in the following natural way:
\begin{itemize}
\item Replace each coloured vertex of $G$ by a clique of order $t$ with every edge within it that colour;
\item Replace each coloured edge of $G$ by a complete bipartite graph $K_{t,t}$ with every edge that colour.
\end{itemize}
Finally, given $t \geq 2$, we call an edge-coloured clique $G'$ \emph{$(r,t)$-unavoidable} if $G'$ is a $t$-blowup of \textit{some} $r$-minimal graph. 
We denote by $\mathcal{F}^r_t$ the set of all $(r,t)$-unavoidable graphs.
Note that putting $r=2$ recovers the definition of $(2,t)$-unavoidable from earlier, 
and also that one can show $G$ being $r$-minimal implies $G$ has at most $2r$ vertices, and so $\mathcal{F}^r_t$ is a finite set (see the argument in~\cite{BLM}).

Bowen, Lamaison and M\"uyesser~\cite{BLM} proved the generalisation of the result of Cutler and Mont\'agh to $r$ colours. 
Not long after, Bowen, Hansberg, Montejano and M\"uyesser~\cite{BHMM} attempted to generalise the result of the first author and Narayanan to $r$ colours, showing that for any $r,t \in \mathbb{N}$, any $r$-edge colouring of a sufficiently large graph on $n$ vertices with $\Omega(n^{2-1/tr^r})$ edges in each colour contains a member of $\mathcal{F}^r_t$.
They conjectured that $\Omega(n^{2-1/t})$ edges in each colour is sufficient. 
Our main result of this section is to confirm this conjecture.

\begin{theorem}\label{thm:unavoidable}
For every $r,t\in \mathbb{N}$, there exists $C(r,t)>0$ such that any $r$-edge colouring of a sufficiently large complete graph on $n$ vertices with the property that any colour appears in at least $C(r,t)n^{2-1/t}$ edges must contain an $(r,t)$-unavoidable graph. 
\end{theorem}


\subsection{Proof of Theorem~\ref{thm:unavoidable}} \hfill\\
We require the following version of the dependent random choice lemma. 
\begin{lemma}[\cite{FS-drc}]\label{lemma:drc}
For all $L,K,t \in \mathbb{N}$, there exists a constant $C$ such that any graph with at least $Cn^{2-1/t}$ edges contains a set $T$ of $L$ vertices, which itself contains a subset $S$ of $K$ vertices in which each subset $X \subseteq S$ with $t$ vertices has a common neighbourhood of size at least $K$. Moreover the neighbourhoods may be chosen so that they are pairwise disjoint from each other and from $T$.  
\end{lemma}


We will also repeatedly apply the bipartite version of  Ramsey's Theorem; since we do not attempt to optimise the bound on the constant $C$ which we will obtain, the following statement suffices for our proof. We write it exactly as we will apply it.
\begin{lemma}\label{lemma:monochr}
For all $r,s,t \in \mathbb{N}$ with $s \leq t$, there exists $s',t'$ with $s' \leq t'$ such that however 
the edges of two cliques $K_{s'}$ and $K_{t'}$ and the edges between them are $r$-coloured, 
there exists a monochromatic $K_s \subseteq K_{s'}$ and a monochromatic $K_t \subseteq K_{t'}$,
with all edges between $K_s$ and $K_t$ monochromatic (where we may have different colours for each of the different monochromatic parts).
\end{lemma}

The idea of the proof is as follows. 
We initially apply the dependent random choice lemma to find large sets $A_i$, $i \in [r]$, such that every subset $T$ of size $t$ of each $A_i$ forms a large monochromatic bipartite graph in colour $i$ with some set $C_i(T)$.
The aim will be to find subsets $D_i \subseteq A_i$, which together with $C_i(D_i)$ form a graph which will then easily be shown to contain an $(r,t)$-unavoidable colour pattern. 
First we apply Lemma~\ref{lemma:monochr} to the pairs $(A_i,A_j)$ so that we have large monochromatic graphs between each pair.
Now we get to the crucial step in our proof which results in the improvement on the bound given in~\cite{BHMM};
we apply Lemma~\ref{lemma:monochr} between $A_j$ and \emph{all} sets $C_i(T)$ with $i<j$. 
Only then do we fix a $t$-set $D_j$ within $A_j$. 
We fix the sets $D_j$, $j=r,r-1,\dots,1$ one at a time, each time applying Lemma~\ref{lemma:monochr} to $C_j(D_j)$ and $A_i$, $i<j$, before going on to fix $D_{j-1}$. 
Finally, we apply Lemma~\ref{lemma:monochr} to each pair $(C_i(D_i),C_j(D_j))$.

The large number of applications of Lemma~\ref{lemma:monochr} results in the $A_i$'s needing to be extremely large, and also increasing in size, i.e. $|A_i|=m_i$, with $1/m_r \ll \dots \ll 1/m_1$.
However, crucially, we can skip the step required in~\cite{BHMM}, to obtain monochromatic subgraphs between the $D_j$'s and $C_i(T)$'s, which is the step which led to the unwanted $r^r$ factor.
Also note that by finding a $t$-blowup of a vertex and edge coloured graph $H$ with a rainbow matching, we got around the issue of having to analyse the structure of $(r,t)$-unavoidable colour patterns.

Note that by $a \ll b$ we mean that there exists a non-decreasing function $f:\mathbb{R}^+ \to \mathbb{R}^+$ such that whichever desired statement we want holds for all $a,b$ satisfying $a \leq f(b)$.

\begin{proof}[Proof of Theorem~\ref{thm:unavoidable}]
Fix beforehand positive constants 
\begin{align}\label{eq:constants}
1/m_r\ll \ldots \ll 1/m_2 \ll 1/m_1\ll 1/r,1/t, \quad 1/C\ll 1/\ell \ll 1/r,1/t.
\end{align}
Suppose that $G$ is an $r$-edge coloured complete graph on $n$ vertices such that every colour appears in at least $Cn^{2-1/t}$ edges. 
We wish to show that $G$ contains an $(r,t)$-unavoidable graph.

First, using Lemma~\ref{lemma:drc}, we may find $r$ pairwise disjoint sets $A_1,\ldots, A_r$ such that 
\begin{enumerate}[(i)]
    \item For every $i\in [r]$, $|A_i|\coloneqq m_i$,
    \item For every $i\in [r]$ and every subset $T\subset A_i$ of size $t$, there is a set $C_i(T)$ such that $(T,C_i(T))$ forms a complete bipartite graph in colour $i$, $|C_i(T)|=\ell$, $C_i(T)$ is vertex disjoint from $\bigcup_{i=1}^{r} A_i$, and the $C_i(T)$ are all pairwise disjoint from each other.
\end{enumerate}
While not stated by Lemma~\ref{lemma:drc} explicitly as written, it is clear that we can guarantee that the sets $A_i$ are pairwise disjoint and also the sets $C_i(T)$ are disjoint from each $A_j$ where $j \not=i$, by using $L$ sufficiently large.

Now we proceed to the main stage of the argument, where we iteratively apply Lemma~\ref{lemma:monochr} to find subsets of the $A_i$'s and subsets of the $C_i(T)$'s which satisfy certain properties.
For ease of notation, we will keep calling these sets $A_i$ and $C_i(T)$. 
We also note that in every step the relations (given by~(\ref{eq:constants})) between the sizes of the new sets are maintained; that is, we make the $m_i$ and $\ell$ sufficiently large so that (\ref{eq:constants}) is maintained after each application of Lemma~\ref{lemma:monochr}. 

\begin{enumerate}
\item For each $i,j \in [r]$ with $i<j$, apply Lemma~\ref{lemma:monochr} to the pair $(A_i,A_j)$. 
\item For each $j \in \{2,\dots,r\}$, for each $i \in [j-1]$, and for each subset $T$ of size $t$ of $A_i$, apply Lemma~\ref{lemma:monochr} to the pair $(A_j,C_i(T))$.
\item For $j$ from $r$ to $1$ (start at $r$ and decrease by $1$ in the iteration), fix $D_j \subseteq A_j$ of size $t$, let $F_j:=C_j(D_j)$, and for each $i \in [j-1]$, apply Lemma~\ref{lemma:monochr} to the pair $(A_i,F_j)$.  
\item For each $i,j \in [r]$ with $i<j$, apply Lemma~\ref{lemma:monochr} to the pair $(F_i,F_j)$. 
\end{enumerate}

As mentioned, after each application of Lemma~\ref{lemma:monochr} we rename everything while maintaining the relations between the sizes of the sets in~(\ref{eq:constants})); for example in (1), when we apply Lemma~\ref{lemma:monochr} to the pair $(A_i,A_j)$, what happens is the following: Ramsey's Theorem guarantees graphs $A_i' \subset A_i$, $A_j' \subset A_j$ such that there is a monochromatic complete bipartite graph between $A_i'$ and $A_j'$, and also $A_i'$ and $A_j'$ are each monochromatic complete graphs. We have $|A_i'|=m_i'$, $|A_j'|=m_j'$ where we have 
$$ 1/m_r  \ll \dots \ll 1/m_{j+1}  \ll 1/m_j' \ll  1/m_{j-1} \ll \dots \ll  1/m_{i+1}\ll  1/m_i' \ll 1/m_{i-1} \ll \dots \ll 1/m_1.$$
Then we relabel so that $A_i:=A_i'$, $A_j:=A_j'$, $m_i:=m_i'$ and $m_j:=m_j'$.
We do similar relabellings for (2), (3) and (4).

Following steps (1) and (2), in addition to properties (i) and (ii) above, 
we have the following:
\begin{enumerate}
\item[(iii)] For every $i, j \in [r]$ with $i<j$, $(A_i,A_j)$ forms a complete monochromatic bipartite graph.
\item[(iv)] For every $i,j \in [r]$ with $i<j$ and every subset $T$ of size $t$ of $A_i$, $(A_j,C_i(T))$ forms a complete monochromatic bipartite graph.
\end{enumerate}

Now we examine the result of step (3). By construction, for each $j \in [r]$, we have
\begin{enumerate}
\item[(v)] By (1), $D_j$ is a subset of $A_j$ of size $t$ which is monochromatic.
\item[(vi)] By (ii), $(D_j,F_j)$ forms a complete monochromatic (in colour $j$) bipartite graph. 
\item[(vii)] By (iv), $(D_j,C_i(T))$ forms a monochromatic complete bipartite graph for every $i<j$ and every $T \subset A_i$ of size $t$. 
In particular, $D_j$ will form a monochromatic complete bipartite graph with $F_i$, once $D_i$ and hence $F_i$ is chosen.
\item[(viii)] By (3), for each $k>j$, $(A_j,F_k)$ and therefore $(D_j,F_k)$ forms a monochromatic complete bipartite graph.
\end{enumerate} 

Finally we examine the result of step (4). We have
\begin{enumerate}
\item[(ix)] For all $i \in [r]$, $F_i$ is monochromatic.
\item[(x)] For every $i, j \in [r]$ with $i<j$, $(F_i,F_j)$ forms a complete monochromatic bipartite graph.
\end{enumerate}

The overall result is that the subgraph $H$ of $G$ induced by $\{D_i,F_i: i \in [r]\}$ forms a $t$-blowup of a vertex and edge coloured clique $J$ on $2r$ vertices where each of the $r$ colours appears at least once. 
Therefore by deleting vertices from $J$ if necessary, we can obtain an $r$-minimal graph.
Hence $G$ contains as a subgraph a member of $\mathcal{F}^r_t$, as required.
\end{proof}

\section{concluding remarks}\label{sec:conc}
We have shown that for any $3$-chromatically connected graph $H$, and any $s\in \mathbb{N}$, there are graphs $G$ for which $G \stackrel{r}{\to} H$ and yet $G[s]\not \stackrel{r}{\to} H[2]$. 
As previously observed, this implies that for these graphs, their blowup Ramsey number depends on the ground graph $G$.
On the other hand, we showed that for forests there is no dependence on the ground graph $G$. 

The first natural problem is to classify all finite graphs for which their blowup Ramsey number does not depend on the ground graph. We conjecture that forests are the only such graphs. 

\begin{conjecture}
Let $H$ be a graph containing a cycle. Then, there exists integers $r,t \geq 2$ such that for every integer $s \geq 2$, there exists a graph $G$ such that $G \stackrel{r}{\to} H$ but $G[s]\not \stackrel{r}{\to} H[t] $. 
\end{conjecture}

It would already be interesting to prove the conjecture for cycles $C_k$ with $k \geq 4$ (which are not $3$-chromatically connected).

Regarding Theorem~\ref{thm:treesmain} it is natural to ask whether the function $f(T,2,t)$ is exponential in $t$, as the result of Fox, Luo and Wigderson suggests. A slightly more careful analysis of our argument shows that we may take $f(T,2,t)=2^{2^{O(t)}}$ but we suspect single exponential might be the truth. 

\begin{problem}
Is $f(T,2,t)=2^{O(t)}$? 
\end{problem}

\section*{Acknowledgements}
The authors would like to thank the two anonymous referees for their careful and helpful reviews.

\bibliographystyle{plain}
\bibliography{tworamsey}

\end{document}